\documentclass[11pt,a4paper]{amsart}
\usepackage{amssymb,amsmath}
\input{amssym.def}
\input{amssym}
\usepackage{graphicx, tikz,caption}
\usepackage{mathtools}
\usepackage{comment}

\newcommand{\N}{{\mathbb N}}

\newcommand{\R}{{\mathbb R}}

\newcommand{\Z}{{\mathbb Z}}

\newcommand\cA{{\mathcal A}}
\newcommand\cB{{\mathcal B}}

\newcommand\cO{{\mathcal O}}

\newcommand\am{\mathfrak a}

\newtheorem{theorem}{Theorem}[section]

\newtheorem{problem}[theorem]{Problem}
\newtheorem{corollary}[theorem]{Corollary}
\newtheorem{lemma}[theorem]{Lemma}
\newtheorem{remark}[theorem]{Remark}
\newtheorem{definition}[theorem]{Definition}
\newtheorem{example}[theorem]{Example}
\newenvironment{proofof}[1]{\noindent {\bf Proof of #1.}}{ \hfill $\Box$\\ }

\title{Rotated Odometers and Actions on Rooted Trees}
\author{Henk Bruin}
\address{Henk Bruin, Faculty of Mathematics, University of Vienna, Oskar-Morgenstern-Platz 1, 1090 Vienna, Austria}
\email{henk.bruin@univie.ac.at}
\author{Olga Lukina}

\address{Olga Lukina, Faculty of Mathematics, University of Vienna, Oskar-Morgenstern-Platz 1, 1090 Vienna, Austria, and Mathematical Institute, Leiden University, P.O. Box 9512,
2300 RA Leiden,
The Netherlands}
\email{o.lukina@math.leidenuniv.nl}
\thanks{This research is supported by the FWF Project P31950-N35}
\date{October 21, 2022}

\begin{document}

\baselineskip=17pt

\maketitle

\begin{abstract}
A rotated odometer is an infinite interval exchange transformation (IET) obtained as a composition of the 
von Neumann-Kakutani map and a finite 
IET of intervals of equal length. In this paper, we consider rotated odometers for which the finite 
IET is of intervals of length $2^{-N}$, for some $N \geq 1$. We show that every such system is 
measurably isomorphic to a $\Z$-action on a rooted tree, and that the unique minimal aperiodic subsystem 
of this action is always measurably isomorphic to the action of the adding machine. 
We discuss the applications of this work to the study of group actions on binary trees.
\end{abstract}

\section{Introduction}\label{sec:intro}

In this paper, we consider infinite interval exchange transformations (IETs) obtained by precomposing 
the von Neumann-Kakutani map of an interval 
with a finite IET of equal length intervals, 
and study the dynamics of such systems. 

Let $\am$ be the von Neumann-Kakutani map, represented on the half-open unit interval $[0,1)$ as
\begin{align}\label{eq-odometer}
\am(x) =  x - (1-3 \cdot 2^{-n}) \qquad  \text{ if } x \in [1-2^{1-n}, 1-2^{-n}),\ n \geq 1.
\end{align}
For $q \in \N$, divide the interval $I = [0,1)$ into $q$ half-open subintervals of length $\frac{1}{q}$. Let $\pi$ be a permutation of $q$ symbols and let $R_\pi$ be the corresponding piecewise continuous map of the subintervals. 
The infinite IET $F_\pi:I \to I$ defined by $F_\pi = \am \circ R_\pi$ is called the \emph{rotated odometer}. 
This generalizes the case when $R_\pi:x \mapsto x+ p/q \mod 1$ is a circle rotation, and we keep the name for
the general case. 

It was shown in \cite{BL2021} that every rotated odometer $(I,F_\pi,\lambda)$ with Lebesgue measure $\lambda$ is measurably isomorphic to the first return map of a flow of rational slope on a certain infinite-type translation surface. 
The translation surfaces in question have interesting properties: they are non-compact surfaces of finite area, 
infinite genus and with a finite number of ends. The closure of such a surface contains a single \emph{wild} singularity and possibly a finite number of cone angle singularities, 
see \cite{DHP,Rbook} for definitions and details about translation surfaces of infinite type. On the other hand, one can consider $(I,F_\pi,\lambda)$ as a perturbation of the von Neumann-Kakutani  system $(I,\am,\lambda)$. 
A natural question is, what dynamical properties of $(I,\am, \lambda)$ are preserved under such perturbation? For the case $q \ne 2^N$, $N \geq 1$, 
this question was partially answered in \cite{BL2021}. 

Let $I_{per}$ the set of periodic points in $I$ and $I_{np} = I \setminus I_{per}$ be the non-periodic points. 
It was shown in \cite{BL2021} that the aperiodic subsystem $(I_{np},F_\pi)$ of the rotated odometer $(I,F_\pi)$ can be embedded into the Bratteli-Vershik system on a suitable Bratteli diagram, which can be constructed using coding partitions.
The ergodic measures and the spectrum of the Koopman operator 
for $(I_{np}, F_\pi)$ can then be studied using the methods developed in the 
literature for stationary Bratteli diagrams, see \cite{BKMS2010,Fogg2002}. 
In \cite{BL2021} we investigated these questions for the case $q \ne 2^N$, $N \geq 1$. 
In particular, it was shown that $(I_{np},F_\pi)$ may be non-minimal with unique minimal set, 
and that it admits at most $q$ invariant ergodic measures
(examples of rotated odometers with $2$ invariant ergodic measures are given too).

In this paper, we consider the case $q = 2^N$, $N \geq 1$, where it is possible to construct a different, simpler Cantor model for the dynamical system of a rotated odometer than in \cite{BL2021}. 
More precisely, we show that the rotated odometer $(I,F_\pi,\lambda)$ is measurably isomorphic to a $\Z$-action on a rooted binary tree, and, using this model, we study the dynamical and ergodic properties of the system. We also discuss the applications of our results to the study of group actions on binary trees.

We now give an overview of the main steps in the procedure which builds a measurable isomorphism between $(I,F_\pi,\lambda)$ and a $\Z$-action on a tree. 

As a first step, we embed $(I,F_\pi)$ into a dynamical system given by a homeomorphism of a Cantor set, that is, there exists a Cantor set $I^*$, a homeomorphism $F_\pi^*:I^* \to I^*$ and an injective map $\iota: I \to I^*$, such that the image $\iota(I)$ is dense in $I^*$ and $\iota \circ F_\pi = F_\pi^* \circ \iota$.  
This procedure has an important difference with an embedding of $(I, F_\pi)$ 
into a compact space $(I^*,F_\pi)$ constructed in \cite{BL2021}.

Indeed, to define the compact space $I^*$ in \cite{BL2021} we employ a technique standard in the study of finite IETs, see for instance \cite{Keane1975}.
Namely, we create gaps in $I$ by doubling points in the orbits of discontinuities of $F_\pi$. 
Periodic points in \cite{BL2021} have half-open neighborhoods where each point is periodic with the same period as $x$, and 
no points in this neighborhood get doubled. Consequently $I^*$ is not totally disconnected. 
However, the closure of $\iota(I_{np})$ is always a Cantor set.

In this paper $I^*$ is constructed by simply doubling every dyadic rational 
$p/2^m$, $m \geq 1$, $0 < p < 2^m$, thus repeating the construction of the middle-third Cantor set, if we think of the middle interval as collapsed to a point. The compact space $I^*$ obtained this way is always totally disconnected. The discontinuity points of $(I,F_\pi)$ are among the doubled points, which implies that $F_\pi$ extends to a homeomorphism $F_\pi^*$ of $I^*$. The embedding $\iota$ is a measurable map with respect to the Lebesgue measure $\lambda$ on $I$ 
and the measure $\mu$ on $I^*$ defined in Section~\ref{subsec-embedding}.

\medskip
We next build a tree model.

\begin{definition}\label{def-binary}
A \emph{rooted binary tree} $T$ consists of the set $V = \bigsqcup_{i \geq 0} V_i$ of vertices and the set $E = \bigsqcup_{i \geq 1}E_i$ of edges, 
which satisfy the following properties for all $i \geq 0$: 
\begin{enumerate}
\item The cardinality $|V_i| = 2^i$.
\item Every vertex in $V_i$ is connected by edges to precisely two vertices in $V_{i+1}$.
\item Every vertex in $V_{i+1}$ is connected by an edge to precisely one vertex in $V_i$.
\end{enumerate}
\end{definition}

We modify the binary tree to obtain a \emph{grafted} binary tree as follows. 

\begin{definition}\label{def-grafted}
For $N \geq 1$, a \emph{grafted} binary tree $T_N$ consists of the set $V = \bigsqcup_{i \geq 0} V_i$ of vertices and the set $E = \bigsqcup_{i \geq 1}E_i$ of edges, such that:
\begin{enumerate}
\item  $|V_0|=1$, $|V_1| = 2^N$ and for $i \geq 2$ we have $|V_i| = 2^{N+i-1}$. 
\item The root $v_0 \in V_0$ is connected by edges to $2^N$ vertices in $V_1$. 
\item For $i \geq 1$, every vertex in $V_i$ is connected by edges to precisely $2$ vertices in $V_{i+1}$, and to a single vertex in $V_{i-1}$.
\end{enumerate} 
\end{definition}

In the notation of Definition~\ref{def-grafted}, we have $T_1 = T$, where $T$ is the binary tree of Definition~\ref{def-binary}.


We introduce a labelling of vertices in $V$. Write $\cA_{k} = \{0,1,\ldots,2^k-1\}$, for $k \geq 1$, 
and consider the tree $T_N$. The root $v_0 \in V_0$ is not labelled, vertices in $V_1$ are labelled by digits in $\cA_{N}$, and for $i \geq 1$, if $v \in V_i$ is labelled by a word $w_1w_2 \cdots w_i$ where $w_1 \in \cA_{N}$ and $w_i \in \cA_1$ for $i \geq 2$, then the two vertices in $V_{i+1}$ connected to $v$ are labelled by $w_1 \cdots w_i 0$ and $w_1 \cdots w_i 1$. 

\begin{definition}
An \emph{infinite path} in the tree $T_N$ is an infinite sequence in the product space
 \begin{align}\label{eq-pathspaceproduct}\partial T_N = \{(w_i) = w_1 w_2 \ldots \mid w_1 \in \cA_{N}, w_i \in \cA_1, i\geq 2\} = \cA_{N} \times \prod_{i \geq 2} \cA_{1,i}, & & \cA_{1,i} = \cA_1 \textrm{ for } i\geq 2. \end{align}
 The space $\partial T_N$ is called the \emph{boundary} of the tree $T_N$.
\end{definition}

Since $N$ is finite and the cardinality of $\cA_1$ is two, $\partial T_N$ is a Cantor set.

\begin{definition}\label{def-auto}
An automorphism $g: T_N \to T_N$ is a map of $T_N$ which restricts to bijective maps on the sets $V$ and $E$ of vertices and edges respectively, and which preserves the structure of the tree. That is, if $v_1 \cdots v_i \in V_i$ is a vertex, then for any vertex $v_1 \cdots v_i w \in V_{i+1}$, where $w \in \{0,1\}$, we have that $g(v_1 \cdots v_i)$ is a subword of
$g (v_1 \cdots v_i w_i)$. In other words, two vertices in $V_i$ and $V_{i+1}$ are joined by an edge if and only if their images under $g$ are joined by an edge.
\end{definition}

We denote by $Aut(T_N)$ the group of automorphisms of $T_N$.  
It is straightforward to see that every automorphism $g \in Aut(T_N)$ induces a homeomorphism of the boundary $\partial T_N$. 

A \emph{cylinder}, or a \emph{cylinder set} $[w_1w_2 \ldots w_i]$ in $T_N$, $i \geq 1$, is the set of all infinite paths starting with the finite sequence $w_1 w_2 \ldots w_i$.
The {\em Bernoulli measure} $\mu_N$ on $\partial T_N$ is the
standard measure in which every cylinder $[w_1w_2\dots w_i]$
has the mass $2^{-N-i+1}$. It is straightforward that $\mu_N$ is preserved under every 
automorphism of $T_N$.

\begin{theorem}\label{thm-main3}
Let $q = 2^N$, let $\pi$ be a permutation on $q$ symbols and let $(I,F_\pi,\lambda)$ be a rotated odometer 
with Lebesgue measure $\lambda$. Then there exists an automorphism $\widetilde F_\pi \in Aut(T_N)$ 
and a measurable isomorphism 
  $${\phi}: (I,F_\pi,\lambda) \to (\partial T_N,\widetilde F_\pi, \mu_N)$$
such that $\widetilde F_\pi \circ \phi = \phi \circ F_\pi$.   
\end{theorem}

Theorem~\ref{thm-main3} is proved in Section~\ref{subsec-treemodel}.

A consequence of Theorem~\ref{thm-main3} is the following description of the dynamics 
of $(I,F_\pi,\lambda)$ in the case $q=2^N$, $N \geq 1$, 
which is more precise than the result of \cite{BL2021}. 

\begin{theorem}\label{thm-main1}
Let $q = 2^N$ for some $N \geq 1$, and let $(I,F_\pi)$ be a rotated odometer.
There exists a decomposition $I = I_{per} \cup I_{np}$ with the following properties:
\begin{enumerate}
\item[(i)] Every point in $I_{per}$ is periodic, the restriction $F_\pi: I_{per} \to I_{per}$ is well-defined and invertible.
\item[(ii)] If $I_{per}$ is non-empty, then $I_{per}$ is a finite union of half-open maximal periodic intervals $[x,y)$, $x,y \in I$. Thus the set of periods of points in $(I,F_\pi)$ is finite.
\item[(iii)] The set $I_{np}$ contains $0$ and $F_\pi: I_{np} \to I_{np}$ is well-defined 
and invertible at every point in $I_{np} \setminus \{ 0 \}$.
\item[(iv)] The aperiodic system $(I_{np}, F_\pi)$ is minimal.
\end{enumerate}
\end{theorem}

The difference with the general case $q \geq 2$ in \cite{BL2021} is that there $I_{per}$ can be an infinite union of half-open intervals, 
while for $q = 2^N$, $I_{per}$ is at most a \emph{finite} union of half-open intervals. It follows that the set of periods which occur in $(I,F_\pi)$ is finite, 
which need not be the case in \cite{BL2021}. Another difference is that for $q = 2^N$ the aperiodic subsystem $(I_{np},F_\pi)$ is always minimal, while this need not hold for $q \ne 2^N$. Theorem~\ref{thm-main1} is proved in Section~\ref{subsec-treemodel}.

Since $I_{per}$ is a finite union of half-open intervals, $I_{np}$ is also a finite union of half-open intervals, 
and its Lebesgue measure $\lambda(I_{np}) > 0$. 
We normalise $\lambda_{np}(U) = \lambda(U)/\lambda(I_{np})$ for every $U \subset I_{np}$. 
The \emph{ dyadic adding machine} $\am: \{0,1\}^\N \to \{0,1\}^\N$ is a well-known example of a minimal $\Z$-action on the space of one-sided infinite sequences of $0$'s and $1$'s. For a finite set $S = \{0,\ldots,r-1\}$, $r \geq 1$, we define the adding machine $\am_S: S \times \{0,1\}^\N \to S \times \{0,1\}^\N $
as the addition of $1$ in $S$ with infinite carry to the right. 
In other words, 
\begin{align}\label{addingS}
\am_S(s,x) = \begin{cases}
                          (s+1,x) & \text{ if } s < r-1,\\[1mm]
                          (0,\am(x)) & \text{ if } s=r-1, \text{ where $\am$ is the dyadic adding machine},
                         \end{cases}
\end{align}
The adding machine $\am_S$ preserves the obvious Bernoulli measure $\mu_S$.

 Since $\{0,1\}^\N$ is homeomorphic to $\partial T_1$, there is a conjugate action on $\partial T_1$ which we also call the adding machine and denote by $a$.
A recursive definition of the adding machine on the boundary $\partial T_N$ of the grafted tree $T_N$ is given in Example~\ref{ex-addmach1}. 

\begin{corollary}\label{cor-addmach}
Let $q = 2^N$ for some $N \geq 1$, and let $(I,F_\pi)$ be a rotated odometer.
The aperiodic system $(I_{np},F_\pi,\lambda_{np})$ is measurably isomorphic to the action of the adding machine on $S \times \{0,1\}^{\N}$, for $|S| \leq 2^N$, with Bernoulli measure $\mu_S$.
\end{corollary}

A rotated odometer need not be conjugate to an automorphism of the binary tree $T$, since  $F_\pi$ may be such 
that the permutation $\pi$ does not respect the structure of the binary tree, see Remark~\ref{rem-period3}. 
Therefore  $T_N$ cannot be substituted by $T_1$ in Theorem~\ref{thm-main3}.

Theorem~\ref{thm-main3} has applications in the study of group actions on binary trees. Infinite IETs and actions of self-similar groups on binary trees are related. For instance, the famous Grigorchuk group was initially defined as a group of infinite IETs 
of the unit interval, see \cite[Section 2]{Grig2011}. Actions of self-similar groups on binary rooted trees are an active topic of research in 
Geometric Group Theory \cite{BN2008,Grig2011,Nekr2005}, and they also have applications in the study of arboreal representations of absolute 
Galois groups of number fields \cite{Jones2013,Lukina2018}. We now present a corollary of Theorem~\ref{thm-main3} for the actions of groups on binary trees. 

To this end, let $T_1 = T$ be the binary tree.
An automorphism $g \in Aut(T)$ is of finite order if $g^m = id$ for some $m \geq 1$. For instance, if $g_i$ interchanges $0$'s and $1$'s in the $i$-th coordinate $w_i$, then $g_i$ has order $2$. Another example of an element of order $2$ is $g_{even}$, which interchanges $0$ and $1$ in $w_i$ for every even $i$, and of course one can construct many more examples. The adding machine \eqref{addingS} is an automorphism of $T$ of infinite order. 

Let $G \subset Aut(T)$ be a profinite group such that $G$ acts transitively on $\partial T$. Given $g \in Aut(T)$, the restriction $g|V_n$ is a permutation of a finite set $V_n$, and so it can be written as a product of cycles. Let $(x_i) = x_1 x_2 \cdots \in \partial T$, then $x_1 \cdots x_n$ is a vertex in $V_n$. Denote by $g_{n,x_1 \cdots x_n}$ the cycle containing $x_1 \cdots x_n$, then one can ask how the sequence of cycles $\{g_{n,x_1 \cdots x_n}\}$ behaves as $n$ increases. 
It is conjectured in \cite{BJ2007}, that when $G$ is a representation of the absolute Galois group of a number field, elements with a certain type of cycle structure are dense in $G$. To the best of our knowledge, this conjecture is solved only in a few cases.

As a rule, given $g \in Aut(T)$, it is not immediate to determine the cycle structure of $g$, except in a few simple cases when $g$ is periodic or when $g$ acts transitively on every level $V_n$, $n \geq 1$. 
The theorem below allows us to determine the cycle structure for compositions of the adding machine and some  periodic elements of $Aut(T)$.

\begin{theorem}\label{theorem-appl}
Let $\mu_1$ be the Bernoulli measure on $\partial T$, and let $\lambda$ be Lebesgue measure on the half-open unit interval $I$. Let $g\in Aut(T)$ be such that there exists $m \geq 1$ such that for every $i > m$ and every sequence $w_1 w_2 \cdots \in \partial T$ the action of $g$ leaves $w_i$ unchanged (which implies that $g$ has finite order). Let $a \in Aut(T)$ be the adding machine. Then the following is true:
\begin{enumerate}
\item For some permutation $\pi$ on $2^m$ intervals, there exists a rotated odometer $(I,F_\pi)$ and an injective measure-preserving map $ \phi: (I,\lambda) \to (\partial T,\mu)$, 
such that $\phi \circ F_\pi = (a \circ g)\circ \phi$.
\item Consequently, $a \circ g$ has infinite order, there is a clopen subset $U \subset \partial T$ such that the restriction $\langle a \circ g \rangle|U$ is minimal, and there is an $n_0 \geq 0$ such that 
every $x \in \partial T \setminus U$ is periodic of period $2^k$ for some $k \leq n_0$.
\end{enumerate}
\end{theorem}

The realization of a tree automorphism as an interval exchange transformation in Theorem~\ref{theorem-appl} relies on the fact that, under the hypotheses of the theorem, $g$ respects the embedding of an interval into the boundary of a tree $T$ in Theorem~\ref{thm-main1}. This means, in particular, that the orbits of points which do not have preimages under $\phi$ consist of points which also do not have preimages under $\phi$. This condition need not hold for a general finite order automorphism of $T$. We discuss this and the possibility of generalizing Theorem~\ref{theorem-appl} to a larger class of tree automorphisms in Remark~\ref{remark-generalize}.

\begin{remark}
{\rm
In the literature, an \emph{odometer} in $Aut(T)$ is sometimes defined as any $h \in Aut(T)$ such that the action of the cyclic group $\langle h \rangle$ is transitive on each $V_n$, $n \geq 1$. Every such $h$ is conjugate to the adding machine in Example~\ref{ex-addmach1} 
by some $g \in Aut(T)$ \cite{Pink13}. We stress that Theorem~\ref{theorem-appl} only holds for the adding machine and need not hold for an odometer $h$. To this end we show in Remark~\ref{remark-notminimal} that it is possible to find $h \in Aut(T)$ such that the action of the cyclic subgroup $\langle h \rangle$ on $\partial T$ is minimal, and a periodic $g \in Aut(T)$, such that the product $h \circ g$ has finite order. 
There exists an infinite IET that is measurably isomorphic to the action of such $\langle h \rangle$ on $\partial T$, but this IET will not be the rotated odometer 
of the form defined at the beginning of the introduction.
}
\end{remark}

We finish with a sample open question motivated by applications to actions on binary trees. 
Consider compositions of the adding machine with a periodic element which does not satisfy the hypotheses of Theorem~\ref{theorem-appl} but which respects the embedding of $I$ in Theorem~\ref{thm-main1}, see Remark~\ref{remark-generalize} for the justification of such an assumption.
It may be possible to solve the following problem by considering a sequence $\{F_{\pi_i}\}_{i \geq 1}$ of rotated odometers, where each $\pi_i$ is a (possibly different) permutation of a finite number of symbols.

\begin{problem}\label{prob-compositions}
Let $g \in Aut(T)$ be periodic such that for any $i \geq 1$ there is $j>i$ and $w_1 \cdots w_j \cdots \in \partial T$ such that $g(w_j) \ne w_j$, and such that $g$ preserves the embedding $\phi$ in Theorem~\ref{thm-main1}. Find a model for the action of the product $a \circ g$,
where $a$ is the adding machine, in terms of rotated odometers.
What are the topological properties of infinite translation surfaces, 
which admit flows whose first return map is measurably isomorphic to such systems?
\end{problem}

The paper is organized as follows. In Section~\ref{sec:general} we develop a tree model for rotated odometers and prove Theorem~\ref{thm-main3}. In Section~\ref{subsec-dynamics} we discuss the dynamics of rotated odometers and prove Theorems~\ref{thm-main1} and~\ref{theorem-appl} and Corollary~\ref{cor-addmach}.

\section{The tree model}\label{sec:general}

In this section we build a tree model for a rotated odometer with $q = 2^N$, and prove Theorem~\ref{thm-main3}.

\subsection{Embedding into a Cantor set}\label{subsec-embedding}

Set $C = \{p 2^{-n} \mid n\geq 1, 0 < p < 2^n\}$; these dyadic rationals are used as cut-points. For each point $x \in C$ we add a double point $x^-$ to $I$, and define
  $$I^* = I \cup \{x^- \mid x \in C\} \cup \{1\}.$$
The subset $I \cup \{1\}$ of $ I^*$ has total order $<$ induced from $\R$. We extend this order to $I^*$ 
by defining $x^- < x$ if $x \in C$, and $y < x^-$ if $y \in I \setminus C$, $x \in C$ and $y < x$. 
Since there are no points between $x^-$ and $x$ in $I^*$, adding $x^-$ to $I$ can be thought of as creating a gap. We give $I^*$ an order topology with open sets
 \begin{align*} \cB = \{(a,b) \mid a,b \in I^*\} \bigsqcup \{[0,b) \mid b \in I^*\} \cup \{(a,1] \mid a \in I^*\}.   \end{align*}
It is straightforward that the sets $\{[x,y^-] \mid x,y \in C\}$ are clopen in this topology. Since $C$ is dense in $I$, every point $z \in I^*$ has a system of decreasing clopen neighborhoods 
 $$C(z,n) = \{[p_n2^{-n},(p_n+1)2^{-n}] \mid n \geq 0, \, 0 \leq p_n < 2^n\}.$$

Recall that a metric $d$ on a space $X$ is an \emph{ultrametric} if it satisfies the following stronger form of the triangle inequality,
$$
d(x,y) = \max\{d(x,z), d(z,y)\} \  \textrm{ for all }x,y,z \in X.
$$
We put an ultrametric on $I^*$ by declaring that
   \begin{align*} d(z_1,z_2) = \frac{1}{2^r}, \quad r = \max\{n \geq 0 \mid C(z_1,n) = C(z_2,n) \}. \end{align*}
Then $I^*$ is a compact totally disconnected perfect metric space, that is, $I^*$ is a Cantor set.

Define a measure $\mu$ on $I^*$ by setting for each clopen set $\{[x,y^-] \mid x,y \in C\}$
  $$\mu([x,y^-]) = y - x,$$
and denote by $\iota: I \to I^*$ the inclusion map. Clearly $\mu(I^*) = 1$. Since $C$ is countable, the following is straightforward.

\begin{lemma}\label{lemma-iota}
The map $\iota: (I,\lambda) \to (I^*,\mu)$ measurable.
\end{lemma}

Denote by $D_0  = \{1 - 2^{-k} \mid k \geq 0\}$ the set of discontinuities of the von Neumann-Kakutani map $\am$, and let $D^+$ and $D^-$ be the sets of forward and backward (whenever defined) orbits of points in $D_0$. Since $\am$ is continuous on the intervals $I_k = [1 - 2^{-(k-1)},1-2^{-k})$, $k \geq 1$, and, moreover, the restriction $\am|I_k$ for each $k \geq 1$ is a translation by $\pm p 2^{-s}$ for some $p,s \in \N$, the set $D_0 \cup D^+ \cup D^-$ of forward and backward orbits of the points of discontinuity of $\am$ is contained in $C$.

We can extend $\am: I \to I$ to a continuous map $\am^*:I^* \to I^*$ by setting $\am^*(x) = \am(x)$ if $x \in I$, and
  $$\am^*(x^-) = \lim_{y \nearrow x} \am(y), \quad \textrm{for all } x \in C \cup \{1\}.$$
Every point $x \in I$ except $0$ has a two-sided orbit, and it follows that $\iota(x)$ has a two-sided orbit in $I^*$. 
 For any sequence $y \nearrow 1$ the sequence of images $\am(y) \searrow 0$, so $\am^*(1) = 0$ and $0$ has a two-sided orbit in $I^*$ under $\am^*$. It follows that $\am^*$ is a homeomorphism. 
 It is immediate that $\am^*  \circ \iota(x) = \iota \circ \am (x)$ for all $x \in I$.

Note that the finite IET $R_\pi: I \to I$ extends in a similar manner to a periodic homeomorphism $R_\pi^*:I^* \to I^*$, which satisfies $R_\pi^* \circ \iota (x) = \iota \circ R_\pi(x)$ for all $x \in I$. Then for the composition $F_\pi^* = \am^* \circ R_\pi^*$ it follows that $F_\pi^* \circ \iota(x) = \iota \circ F_\pi(x)$ for all $x \in I$, and thus $(I,F_\pi,\lambda)$ is measurably isomorphic to $(I^*,F_\pi^*,\mu)$ via the embedding $\iota$. 

\subsection{Actions on trees}

The binary tree and the grafted binary trees were defined in Definitions~\ref{def-binary} and~\ref{def-grafted}, and automorphisms of trees were defined in Definition~\ref{def-auto}. We now introduce a description of elements in $Aut(T_N)$ convenient for computations. This approach is a slight modification of the one routinely used in Geometric Group Theory to study actions on binary trees,
see for instance \cite{Nekr2005}. The purpose of this modification is to take into account the fact that in the grafted binary tree the vertex set $V_1$ has more than $2$ vertices. 

Let $T = T_1$ be the binary tree with the labelling of vertices by finite words in $\cA_1$ as defined in the Introduction. Let $w = w_1 \cdots w_k \in \prod_{i=1}^k\cA_1$ and denote by $T(w)$ the subtree of $T$ consisting of all paths starting with the finite word $w$. All such paths pass through the vertex in $V_k$ labelled by $w$. Then there is an isomorphism of trees 
  \begin{align}\label{eq-binaryhomeo}
  \kappa_w: T(w) \to T, \qquad  w_1 \cdots w_k v_{k+1} \cdots \mapsto v_{k+1} \cdots, \, v_i \in \{0,1\} \textrm{ for }i > k.\end{align}
For every $g \in Aut(T)$, the restriction $g|V_n$ is a permutation of a set of $2^n$ elements. 

\begin{definition}\label{def-section}
Given an automorphism $g \in Aut(T)$, and a finite word $w$, we define a \emph{section} at $w$ by
  \begin{align}\label{eq-section}g_w =  \kappa_{g(w)} \circ g \circ \kappa_w^{-1} \in Aut(T).\end{align}
 \end{definition}
Let $g|V_n = \tau$. Then we can write $g$ as a composition (we compose the maps on the left)
  \begin{align}\label{eq-recur}
  g = (g_{\tau^{-1}(0^n)},g_{\tau^{-1}(0^{n-1}1)}, \ldots, g_{\tau^{-1}(1^n)}) \tau, 
  \end{align}
where $g_{\tau^{-1}(w)}$ are sections, for finite words $w$ of $n$ letters.
Equation \eqref{eq-recur} means that to compute $g$, we first apply $\tau$ on $V_n$, 
and then we apply a section $g_{\tau^{-1}(w)}$ to the subtree $T(w)$, for all $w \in V_n$.
  
\begin{example}\label{ex-addmach1}
{\rm
Using sections, we can write automorphisms of $T$ recursively. Recall that a generator of the adding machine action on a Cantor space $\{0,1\}^{\N}$ is given by
  \begin{align}\label{eq-addingmach}a (w_1 w_2 \cdots) = \left\{ \begin{array}{ll}  (w_1+1)\, w_2 \cdots & \textrm{ if }w_1 = 0, \\  0 \, 0 \cdots 0 \, (w_{k+1}+1) \, w_{k+2} \cdots &\textrm{ if } w_i= 1 \textrm{ for } 1 \leq i \leq k,  \, w_{k+1} =0, \\ 0 \, 0 \cdots  & \textrm{ if }w_k = 1 \textrm{ for all }k \geq 1.  \end{array}\right. \end{align}
Recall that for the binary tree $T$ we have $\partial T \cong \{0,1\}^\N$.   
Let $\sigma$ be the non-trivial permutation of $\cA_1$. Then using \eqref{eq-recur} we can write
  $$a = (a,1)\sigma,$$
where $1$ is the identity map in $Aut(T)$. 
Here $\sigma$ performs the addition of $1$ modulo two in the first entry of the sequence, interchanging $0$ and $1$, while $(a,1)$ implements the recursive procedure of infinite carry to the right. For example, if $w = 10^\infty$, then applying $\sigma$ to $w$ interchanges $1$ to $0$ in the first component, so $\sigma(w) = 0^\infty$, and we must compute $(a,1)(0^\infty)$ next. The sequence $0^\infty$ belongs to the subtree $T(0)$ which means that we must apply the section $a_0=a$ to $0^\infty$. That is, we apply $a$ to $0^\infty$ starting from the second entry. Since $a|V_1 = \sigma$, we must interchange $0$ and $1$ in the second entry, obtaining $010^\infty \in T(01)$. We have for the sections $a_1 = 1$, then also $a_{01} = 1$, and the computation stops with the result $a(10^\infty) = 010^\infty$.  
}
\end{example}  

Using \eqref{eq-recur} we can compute the compositions of elements in $Aut(T)$.  The following statement is obtained by a straightforward computation.

\begin{lemma}\label{lemma-product}
Let $g,h \in Aut(T)$, and suppose $g = (g_0,\ldots,g_{2^n-1}) \tau $ and $h = (h_0,\ldots,h_{2^n-1}) \nu $, where $\tau, \nu$ are permutations of $2^n$ symbols and $g_i,h_i \in Aut(T)$ for $0 \leq i < 2^n$. Then
  \begin{align}\label{eq-compose} gh = (g_0,\ldots,g_{2^n-1}) \tau (h_0,\ldots,h_{2^n-1}) \nu = (g_0 h_{\tau^{-1}(0)},\ldots,g_{2^n-1}h_{\tau^{-1}(2^n-1)}) \tau \nu.\end{align}
  \end{lemma}
  
Now let $T_N$ be the grafted binary tree. Similarly to \eqref{eq-binaryhomeo}, for any $w \in V_k$, $k \geq 1$ we define a map
  \begin{align}\label{eq-kappaw}\kappa_w: T_N(w) \to T_1, \qquad w_1 \cdots w_k v_{k+1} \cdots \mapsto v_{k+1} \cdots.\end{align}
The difference with \eqref{eq-binaryhomeo} is that the range of $\kappa_w$ is not the grafted tree $T_N$ but the binary tree $T$. A section $g_w$ of the grafted tree $T_N$ at $w$ is defined by \eqref{eq-section} with $\kappa_w$ given by \eqref{eq-kappaw}. 
Again, the difference with the setting of the binary tree is that for the grafted binary tree $T_N$ sections are elements of $Aut(T)$ and not of $Aut(T_N)$.

\begin{lemma}\label{lemma-1}
Given an automorphism $g \in Aut(T)$ of the binary tree $T$, there is always an automorphism $\widehat g \in Aut(T_N)$ of the grafted tree $T_N$, such that the induced homeomorphisms on the boundaries of the corresponding trees are conjugate. 
\end{lemma}

\begin{proofof}{Lemma \ref{lemma-1}}
Vertices in the vertex level set $V_N$ of $T$ are labelled by words of length $N$ in the alphabet $\cA_1$. Define the map
  $$\kappa_N: \cA_1^N \to \cA_{N},  \qquad w_1 \cdots w_N \mapsto \sum_{i=1}^N 2^{N-i} w_i.$$
Using the identification \eqref{eq-pathspaceproduct} of the path spaces $\partial T$ and $\partial T_N$ with products of finite sets  we obtain a homeomorphism
  $$\kappa_\infty: \partial T \to \partial T_N, \qquad w_1 w_2 \cdots \mapsto \kappa_N(w_1\cdots w_N)w_{N+1} \cdots.$$
It follows that the map $\widetilde g = \kappa_\infty \circ g \circ \kappa_\infty^{-1} : \partial T_N \to \partial T_N$ is a homeomorphism. Moreover, by construction if two paths $(w_i),(v_i) \in \partial T$ coincide up to level $m \geq N$, then their images under $\kappa_\infty$ coincide up to level $m-N$, so every subtree $T(w)$ for $w \in V_N$ is mapped isomorphically onto a subtree $T_N(\kappa_N(w))$. It follows that $\widetilde g $ defines an automorphism $\widehat g$ of $T_N$.
\end{proofof}

Given a recursive definition of $g \in Aut(T)$ as in \eqref{eq-recur}, we can obtain a recursive definition of $\widehat g \in Aut(T_N)$. Indeed,  let $g|V_N = \tau$ be a permutation of $V_N$ induced by $g$. Then $\tau_N = \kappa_N \circ \tau \circ \kappa_N^{-1}$ is a permutation of the level set $V^N_1$ of $T_N$, and if $g = (g_0,\ldots,g_{2^N-1})\tau$, then $\widehat g =(g_0,\ldots,g_{2^N-1})\tau_N $.

\begin{example}\label{ex-addmach2}
{\rm
Let $a = (a,1)\sigma$ be the standard adding machine as in Example~\ref{ex-addmach1}, and let $N \geq 2$. We can compute that
  \begin{align}\label{eq-tauodometer}\tau_N = \kappa_N \circ (a|V_N) \circ \kappa_N^{-1} = (0, \, 2^{N-1}, \, 2^{N-2}, \, 2^{N-1} + 2^{N-2}, \, \ldots  ,2^N-1),\end{align}
and $\widehat a = (a,1,\ldots,1)\tau_N \in Aut(T_N)$.

More generally, given a finite set $S \subset V_N$, we can consider a subtree $T_S = \bigcup_{w \in S} T_N(w) \subset T_N$. Let $\eta$ be a transitive permutation of $S$, and consider the map $a_S = (a,1,\ldots,1)\eta$ on $T_S$. Then $a_S$ is transitive on $V_n \cap T_S$, for any $n \geq 1$, so $a_S$ is the adding machine on $T_S$.
}
\end{example}

We note that, given $h \in Aut(T_N)$, the composition $\kappa_\infty^{-1} \circ h \circ \kappa_\infty$ need not define an automorphism of $T$. Indeed, let $N=2$, so $T_N$ has $4$ vertices at the first level, and let $h|V_2 = \tau_2 = (012)$, so the vertex $3$ is fixed. We have $\kappa_2^{-1}(3) = 11 \in V_2$ and $\kappa_2^{-1}(2) = 10 \in V_2$. At the same time
  $$ \kappa_2^{-1} \circ \tau_2 \circ \kappa_2(10) = \kappa_2^{-1}(\tau_2(2)) = \kappa_2^{-1}(0) = 00.$$
Thus $\kappa_\infty^{-1} \circ h \circ \kappa_\infty$ maps paths starting with $1$ in $\partial T$ to paths starting with either $1$ or $0$ depending on the second symbol in the sequence. 
This means that $ \kappa_2^{-1} \circ \tau_2 \circ \kappa_2$ is incompatible with the structure of the binary tree $T$, and so $\kappa_\infty^{-1} \circ h \circ \kappa_\infty$ does not define an automorphism of $T$.

\subsection{Tree models for rotated odometers}\label{subsec-treemodel} 
In this section we prove Theorem~\ref{thm-main3}.

\begin{proofof}{Theorem~\ref{thm-main3}}
Recall that $\pi$ is a permutation of $2^N$ symbols, and $\iota: (I,\lambda) \to (I^*,\mu)$ is a measurable embedding into a Cantor set. Write $x_{n,p} = p2^{-n}$ for points in $C$, and $x_{n,p}^-$ for the corresponding double points in $I^*$. 
For each $n \geq 1$, set $x_{n,2^n}^- = 1$. 

Note that for any $n \geq 0$ we have
  $$I^*  = \bigcup \{[x_{n,p}, x_{n,p+1}^-] \mid 0 \leq p < 2^n-1 \}. $$

Consider the grafted tree $T_N$, and recall that $|V_1| = 2^N$. We are going to construct a homeomorphism $\widetilde \phi: I^* \to \partial T_N $ inductively as follows.

Define $\widetilde \phi_1: I^* \to \cA_{N}$ by setting
  $$\widetilde \phi_1(z)  = p \quad \textrm{ if and only if } \quad z \in [x_{N,p}, x_{N,p+1}^-].$$
For $n \geq 2$, there is a unique $0 \leq m < 2^{n}-1$ such that $z \in [x_{n,m}, x_{n,m+1}^-]$. Set 
 $$w_n = \widetilde \phi_{n}(z) = m \mod 2.$$ Then define
  $$\widetilde \phi_\infty: I^* \to \partial T_N, \qquad z \mapsto (\widetilde \phi_1(z), \widetilde \phi_2(z), \ldots).$$
This mapping is bijective, since every point in $I^*$ has a system of clopen neighborhoods of the form $\{[x_{n,p},x_{n,p+1}^-] \mid n \geq 0\}$, and every clopen neighborhood $[x_{n,p}, x_{n,p+1}^-]$ is non-empty. The mapping $\widetilde \phi_\infty$ is clearly continuous and so it is a homeomorphism.   Note that by construction the inclusions of clopen sets in $I^*$ correspond to vertices in $T_N$ joined by finite paths. 

The measure $\mu$ assigns equal weight to each interval $\{[x_{n,p},x_{n,p+1}^-] \mid 0 \leq p <2^n\}$ in the partition of $I^*$, and $\mu(I^*) = 1$. By construction each $[x_{n,p},x_{n,p+1}^-] $ is mapped onto a unique vertex in $V_{n-N+1}$. The Bernoulli measure $\mu_N$ assigns equal weight to every set $\partial T_N(w)$, where $w \in V_{n-N+1}$, and $\mu_N(\partial T) = 1$. It follows that $\widetilde \phi_\infty$ is measure-preserving.

Every map of $I^*$ which for all $n \geq 1$ induces a permutation of clopen sets $\{[x_{n,p}, x_{n,p+1}] \mid 0 \leq p < 2^n-1\}$, induces a family of permutations of the vertex level sets $V_{n-N+1}$, $n \geq 1$ of $T_N$. Since paths in $T_N$ correspond to inclusions of clopen sets in $I^*$, such permutations are compatible with the structure of the tree $T_N$ and induce an automorphism of $T_N$. We note that the maps $\am^*: I^* \to I^*$ and $R_\pi^*:I^* \to I^*$ described in Section~\ref{subsec-embedding} satisfy this condition. Therefore, the composition $F_\pi^* = \am^* \circ R_\pi^*: I^* \to I^*$ induces an automorphism of $T_N$. The proof of Theorem~\ref{thm-main3} is completed by composing $\phi = \widetilde \phi_\infty \circ \iota: I \to \partial T_N$ with the measurable isomorphism $\iota: (I,F_\pi,\lambda) \to (I^*,F_\pi^*,\mu)$.
\end{proofof}

\begin{remark}\label{remark-addedpoints}
{\rm
Consider the set of added points $\{x_{m,p}^- \mid x_{m,p} \in C\}$. Suppose $x_{m,p} = p2^{-m}$ is an irreducible fraction, 
that is, $p$ is odd. Then for $n > m$ we have that $\widetilde \phi_n(x_{m,p}^-) = 1$ since in that case $x_{m,p}^-$ 
corresponds to a right endpoint of a clopen interval in the partition $\{[x_{n,r},x_{n,r+1}^-] \mid 0 \leq r < 2^n\}$, 
and it is always contained in the second interval of the subdivision of $[x_{n,r},x_{n,r+1}^-]$ into two intervals. 
Then the image of $x_{m,p}^-$ in $\partial T_N$ is a sequence which is eventually constant with entries equal to $1$.
}
\end{remark}

\section{Dynamics of rotated odometers} \label{subsec-dynamics}

Using the tree model obtained in Theorem~\ref{thm-main3} we study the dynamics of rotated odometers and prove Theorems~\ref{thm-main1} and~\ref{theorem-appl} and Corollary~\ref{cor-addmach}.

\subsection{Periodic and non-periodic points}\label{subsec-periodic}
For the von Neumann-Kakutani map $\am^*:I^* \to I^*$ denote by $A = \widetilde\phi_\infty \circ \am^* \circ \widetilde\phi_\infty^{-1}: \partial T_N \to \partial T_N$ 
the induced map of the binary tree $T$. We want to describe $A$ using the recursive formula \eqref{eq-recur}. 

\begin{proofof}{Theorem~\ref{thm-main1} and Corollary~\ref{cor-addmach}}
In what follows $n \geq N$.
Let $L_n = [0, 2^{-n})$ and $M_n = [1-2^{-n},1)$, so that $L_n$ is the first and $M_n$ is the last set of the partition of $I$ into $2^n$ sets of equal lengths. Then $\iota(L_n) \subset [x_{n,0},x_{n,1}^-] \subset I^*$ and $\iota(M_n) \subset [x_{n,2^{n}-1},1] \subset I^*$. The definition of the von Neumann-Kakutani map in \eqref{eq-odometer} implies that $\am(x) \in L_n$ if and only if $x \in M_n$, and for any $[x_{n,p}, x_{n,p+1})$ except $M_n$ the restriction $\am|[x_{n,p}, x_{n,p+1})$ is a translation. Thus it preserves the order $\leq $ on the points in $[x_{n,p}, x_{n,p+1})$ induced from $\R$.  
The relation $\leq$ is not preserved by the restriction $\am: M_n \to L_n$, where the order of two halves of $M_n$ is interchanged, and the intervals inside the image of the second half of $M_n$ are further interchanged. The second half of $M_n$ is the set $M_{n+1}$, and we have $\am(M_{n+1}) = L_{n+1}$. Thus the restriction of $\am$ to the set of intervals $\{[x_{n,p}, x_{n,p+1}) \mid 0 \leq p < 2^n\}$, and therefore of $\am^*$ to the set of intervals $\{[x_{n,p}, x_{n,p+1}^-] \mid 0 \leq p < 2^n\}$, defines a permutation of $2^n$ symbols, which is transitive since $\am$ is minimal on $I$. 

It follows that $A|V_{n-N+1}$ is a transitive permutation of $V_{n-N+1}$. Since further permutations of subintervals, which do not respect the order $<$, only happen for the interval, mapped onto $L_{n-N+1}$, for any $w \ne 0^{n-N+1} \in V_{n-N+1}$, 
the section $A_w \in Aut(T)$ is the identity map. The restriction of the section $A_{0^{n-N+1}}$ to $V_{n-N+2}$ is a non-trivial permutation of two symbols, since $\am$ permutes two subintervals of $L_{n-N+1}$. For $n = N$, we have $A|V_{N-N+1} = A|V_1= \tau_N$, for $\tau_N$ given by \eqref{eq-tauodometer}, and so $A = (a,1,\ldots,1)\tau_N$, where $a = (a,1) \sigma$ is described in Example~\ref{ex-addmach1}.

Similarly, given a permutation $\pi$ of $2^N$ symbols, and the corresponding finite IET $R_\pi: I \to I$, we deduce that the induced map $R = \widetilde\phi_\infty \circ R_\pi^* \circ \widetilde\phi_\infty^{-1}$ is given by $R = (1,1,\ldots,1)\pi$, with $R|V_1 = \pi$.

Now using the law for composition of tree automorphisms \eqref{eq-compose} we can easily understand the dynamics of the system $(\partial T_N, A \circ R)$. In particular,
  $$A \circ R = (a,1,\ldots,1)\tau \pi,$$
which leads to the following conclusions:
\begin{enumerate} 
\item[(i)] Consider the decomposition of $\tau \pi$ into cycles, and suppose $c$ is a cycle containing $0$. Let $\cO \subset \cA_{N}$ be the set of symbols in $c$. Then 
  $$\partial T_N(\cO) = \bigcup \{\partial T(s) \mid s \in \cO \}$$ 
is a clopen subset of $\partial T$ and the restriction of $A \circ R$ to this set satisfies
  $$A\circ R|\partial T_N(\cO) = (a,1,\ldots,1)c,$$ 
which shows that this system is the addition of $1$ in the first component with infinite carry to the right, and so it is minimal. Set $S=\{0,\ldots,|c|-1\}$, then $A \circ R| \partial T_N(\cO)$ is isomorphic to the adding machine on $S \times \{0,1\}^\N$ defined in Example~\ref{ex-addmach2}, and Corollary~\ref{cor-addmach} follows. Here $|c|$ denotes the length of the cycle $c$. In particular, the system $(\partial T_N, A \circ R)$ is minimal if and only if $\tau\pi$ is a transitive permutation.
\item[(ii)] In the cycle decomposition of $\tau \pi$, let $c'$ be a cycle not containing $0$, and let $\cO' \subset \cA_{N}$ be the set of symbols in $c'$. Then $\partial T_N(\cO') = \bigcup \{ \partial T(s) \mid s \in \cO' \} $ is a clopen subset of $\partial T$, and we have
 $$A\circ R|\partial T_N(\cO') = (1,1,\ldots,1)c'.$$ 
Thus every point in $T_N(\cO')$ has period $|c'|$. 
\item[(iii)] Since $\tau \pi$ contains a finite number of cycles, the set of periods of periodic points in $(\partial T_N, A \circ R)$, and so in $(I,F_\pi)$, is finite. Also, it follows that a point $x \in \partial T_N$ is periodic if and only if $x \in \partial T_N(\cO')$ for some cycle $c'$ not containing $0$. There is at most a finite number of such cycles $c'$ in $\tau \pi$, and so there is a finite number of half-open intervals in $I$ whose image under the inclusion map $\phi = \widetilde \phi_\infty \circ \iota$ is contained in $\partial T_N\setminus \partial T_N(\cO)$. It follows that the set of periodic points in $(I,F_\pi)$ is at most a finite union of half-open intervals.
\end{enumerate}  
These prove Theorem~\ref{thm-main1} and Corollary~\ref{cor-addmach}.
\end{proofof}

\begin{remark}\label{rem-period3}
{\rm We note that the periods of points in $(I,F_\pi)$ need not be powers of $2$. Let $N = 2$, then $A = (a,1,1,1)(0213)$. Let $\pi = (03)$. Then
  $$A \circ R = (a,1,1,1)(0213)(03) = (a,1,1,1)(0)(321),$$
so the orbit of every infinite sequence in $\partial T_2$ starting with $1$, $2$ or $3$ is periodic with period $3$.

It follows from this example that there exist rotated odometers whose action is not measurably isomorphic 
to the action of an automorphism of the binary tree $T$. 
Indeed, if $g \in Aut(T)$ is an automorphism and $x \in \partial T$ is a periodic point, then the period of $x$ is a power of $2$. 
To see this, consider $x = (x_0,x_1,\ldots)$, and let $r_k$ be the period of $x_k$ in $V_k$. 
Then the period $x_{k+1}$ in $V_{k+1}$ is either $r_k$ or $2r_k$. 
Since the period $r_1$ of $x_1$ in $V_1$ is either $1$ or $2$, the statement follows.
}
\end{remark}

\subsection{Applications}
We prove Theorem~\ref{theorem-appl}.

\begin{proofof}{Theorem~\ref{theorem-appl}}
Suppose that $g \in Aut(T)$ is of finite order such that there is $m \geq 1$ such that for every $i >m$ the action of $g$ leaves $w_i$ unchanged.
We need to show that there exists a rotated odometer $(I,F_\pi)$ for some permutation $\pi$ on $2^m$ intervals, such that $(I,F_\pi,\lambda)$ is measurably isomorphic to $(\partial T, a \circ g , \mu_1)$, 
where $\lambda$ is Lebesgue measure, $\mu_1$ is the Bernoulli measure on the binary tree $T$ and $a$ is the adding machine described in Example~\ref{ex-addmach1}.

Consider the partition of $\partial T$ into clopen sets $\partial T(w)$, where $w= w_1 \cdots w_n$. Also, consider a partition of $I^*$ into subintervals $\{[x_{n,p},x_{n,p+1}) \mid 0 \leq p < 2^n\}$. By construction every such subinterval is mapped under $\phi = \widetilde \phi_\infty \circ \iota$ into a distinct clopen set $\partial T(w)$, and $\phi$ is injective on $I$. 
Define 
$$\widetilde{g}:I \to I, \qquad x \mapsto (\widetilde \phi_\infty \circ \iota)^{-1} \circ g \circ (\widetilde \phi_\infty \circ \iota)(x).$$
 The map $\widetilde{g}$ is well-defined. Indeed, by Remark~\ref{remark-addedpoints} the points in $I^*$ which do not  have preimages in $I$ under $\iota$ are mapped into sequences which are eventually constant with entries equal to $1$. Since $g$ does not change $w_i$ for $i \geq m$, $w \in \partial T$ is eventually a sequence of $1$'s if and only if $g(w)$ is eventually a sequence of $1$'s. Therefore, the map $\phi$ is invertible at $g \circ \phi(x)$. Since $g$ does not change $w_i$ for $i \geq m$, $\widetilde g$ preserves the order of points in the sets $\{[x_{n,p},x_{n,p+1}) \mid 0 \leq p < 2^n\}$, and it follows that the restriction of $\widetilde{g}$ to every interval $\{[x_{m,p},x_{m,p+1}) \mid 0 \leq p <2^m\}$ is a translation. We conclude that $\widetilde g: I \to I$ is a finite IET.
 
It is proved in Section~\ref{subsec-periodic} that the von Neumann-Kakutani map $(I, \am,\lambda)$ is measurably isomorphic to $(\partial T, a, \mu_1)$, where $a = (a,1)\sigma$ is the standard adding machine. Set $F_\pi = \am \circ \widetilde g $, then $(I, F_\pi, \lambda)$ is measurably isomorphic to $(\partial T, a \circ g, \mu_1)$. The second statement of Theorem~\ref{theorem-appl} follows from Theorem~\ref{thm-main1}.
\end{proofof}

\begin{remark}\label{remark-notminimal}
{\rm
In the literature a transformation $g$ such that the cyclic group $\langle h \rangle$ acts transitively on every level $V_n$, $n \geq 1$, of the tree $T$,  is sometimes called an \emph{odometer}. Every odometer is conjugate to the adding machine  in Example~\ref{ex-addmach1} 
by a tree automorphism  \cite{Pink13}. We note that Theorem~\ref{theorem-appl} only holds for the adding machine, but need not hold for its conjugates. 

Indeed, define $a_1 = \sigma$ and $a_2 = (a_1,a_2)$ using the recursive notation, then $\langle a_1,a_2\rangle$ is the dihedral group. Both $a_1$ and $a_2$ have order two, and $h = a_1a_2$ generates an infinite cyclic group whose action on every $V_n$, $n \geq 1$, is transitive. The element $h$ is conjugate to $a = (a,1)\sigma$ but it is not equal to $a$. Recall that $\sigma$ acts on $(w_1,w_2,\ldots) \in \partial T$ by interchanging $0$ and $1$ in the first entry, and keeps the remaining entries fixed. We compute that $h \sigma = \sigma (a_1,a_2) \sigma  = (a_2,a_1)$ has order two.

}
\end{remark}

\begin{remark}\label{remark-generalize}
{\rm
We have seen in the proof of Theorem~\ref{theorem-appl} that, under its hypotheses, a finite order element $g \in Aut(T)$ respects the embedding of $I$ into $\partial T$. More precisely, by Remark~\ref{remark-addedpoints} points which do not have preimages under this embedding correspond to sequences which are eventually constant with entries equal to $1$, and if $g \in Aut(T)$ satisfies the hypotheses of Theorem~\ref{theorem-appl}, then it preserves the set of such sequences.  Suppose $g \in Aut(T)$ does not satisfy the hypotheses of Theorem~\ref{theorem-appl}, that is, for any $n \geq 1$ there is $(w_i) \in \partial T$ and $m_n \geq n$ such that $g(w_{m_n}) \ne w_{m_n}$. Then the action of $g$ on $\partial T$ may or may not respect the embedding of $I$.
If such $g \in Aut(T)$ respects the embedding of $I$ into $\partial T$ ($h$ in Remark~\ref{remark-notminimal} is an example), then $g$ induces an IET of \emph{infinite} number of intervals. At the moment we do not have a unified way of describing the dynamics of a composition of such an IET with the von Neumann-Kakutani map, and we pose this as an open question in Problem~\ref{prob-compositions}. An example of an element which does not respect the embedding is, for instance, $g_{even}$ given for any $(w_i) \in \partial T$ by
  $$g_{even}(w_{2n}) = w_{2n}+1 \mod 2,  \quad g_{even}(w_{2n+1}) = w_{2n+1}, \quad n \geq 1.$$
It is not clear whether such elements induce IETs on the interval $I$, even if one discards the measure $0$ set of orbits in $\partial T$ which do not have preimages under the embedding, and/or allows reflections of subintervals.}
\end{remark}

\end{document}